\documentclass{amsart}
\usepackage{verbatim}
\usepackage[english]{babel}
\usepackage{leftidx}
\usepackage[numbers]{natbib}
\usepackage{bigints}
\usepackage{empheq}
\usepackage{array}
\usepackage{graphicx}
\usepackage{enumerate}
\usepackage{lipsum}
\usepackage{amsthm}
\usepackage{needspace}
\usepackage[makeroom]{cancel}
\usepackage{enumitem}
\newenvironment{myproof}[2] {\paragraph{\textbf{Proof of {#1} {#2} :}}}{\hfill$\square$}
\newtheorem{theorem}{Theorem}[section]
\newtheorem{proposition}[theorem]{Proposition}

\newtheorem{question-non}[]{}

\newtheorem{cor}[theorem]{Corollary}
\newtheorem{definition}[theorem]{Definition}
\newtheorem{observation}[theorem]{Remark}
\newtheorem{example}[theorem]{Example}

\title{On Warped Product Gradient Ricci-Harmonic Soliton}
\author{Batista, E. $^{1}$}
\address{$^{1}$ Instituto Federal do Tocantins, Rodovia TO 040 - Km 349 Loteamento Rio Palmeira, 77300-000, Dian\'opolis, TO Brazil.}
\email{elismardb@gmail.com $^{1}$}

\author{Adriano, L. $^{2}$}
\address{$^{2}$ Universidade Federal de Goi\'as, IME, 131, 74001-970, Goi\^ania, GO, Brazil.}
\email{levi@ufg.br $^{2}$}

\author{Tokura, W. $^{3}$}
\address{$^{4}$ Universidade Federal de Goi\'as, IME, 131, 74001-970, Goi\^ania, GO, Brazil.}
\email{williamisaotokura@hotmail.com $^{3}$}



\date{March 17, 2019}

\begin{document} 
\begin{abstract}In this paper we study gradient Ricci-Harmonic soliton with structure of warped product manifold. We obtain some triviality results for the potential function, warping function and the harmonic map which reaches maximum or minimum. In order to obtain nontrivial examples of warped product gradient Ricci-harmonic soliton, we consider the base and fiber conformal to a semi-Euclidean space which is invariant under the action of a translation group of co-dimension one. This approach provide infinitely many geodesically complete examples in the semi-Riemannian context, which is not contemplated in the Riemannian case by the Theorem 1.2 in \cite{ZHU2017882}.
\end{abstract}
\maketitle
\section{Introduction and main results}

Let $(M^n,g)$ and $(N^m,g')$ be two semi-Riemannian manifolds, $u:M\rightarrow N$ a differential map. If there is a nonnegative constant $\theta$, a smooth function $h : M\rightarrow \mathbb{R}$ and $\lambda \in \mathbb{R}$ such that
\begin{equation}\label{definitiongrhs}
\begin{split}
Ric_{g}+ Hess h -\theta \nabla u \otimes \nabla u&=\lambda g\\
\tau_{g}u-g(\nabla u, \nabla h)&=0,
\end{split}
\end{equation}
where $Ric_{g(t)}$ is the Ricci curvature of the $(M,g(t))$ and $\tau_{g(t)}u=trac_{g}(\nabla d u)$  is the tension field of $u,$ then $((M, g), (N, g'), u, h, \lambda)$ is called a gradient Ricci-Harmonic soliton. For convenience we will call it GRHS. We have that GRHS is shrinking, steady or expanding for $\lambda>0, \lambda=0$ or $\lambda< 0$, respectively.
If $h$ occurs as a constant the soliton is called \textit{trivial}. The function $u$ is called harmonic map, and for $u=constant$ the equation \eqref{definitiongrhs} defines a gradient Ricci soliton metric, which was very important in the Thurston's Geometrization Conjecture studied by Perelman \cite{perelman2002}. For more details on gradient Ricci soliton see \cite{Fernandezlopes2008, MR1207538, MR2480290, MR2507581}. If $(N, g')=(\mathbb{R}, dr^2)$ and $u$, $h$ are constants, then (\ref{definitiongrhs}) defines a Einstein metric, and if $h$ is constant and $u$ is harmonic, then (\ref{definitiongrhs}) defines a harmonic-Einstein metric, which is a natural generalization of Einstein metrics.

After the introduction of the GRHS, much efforts have been devoted to understanding its geometry. For instance, in \cite{ZHU2017882} the authors obtain that any shrinking or steady GRHS is a gradient Ricci soliton if the sectional curvature of $N$ is bounded from above by a constant. On the other hand, in \cite{MR3406677} the authors provide triviality results for the potential function and harmonic map by means of compactness of $M$. 

In recent works on soliton, the notion of warped product introduced in (\cite{bishop1969manifolds}), has become an important tool in the construction of manifolds with certain geometric properties \cite{bang2017differential, dobarro1987scalar, kim2013warped}.

\begin{definition}\label{definitionwarpedproduct}
Let $M=(B^n,g_{B}),(F^m,g_{F})$ two semi-Riemannian manifolds and $f>0$ on $B$. The product manifold $M=B\times F$ furnished with the metric tensor
\begin{equation}\label{metrictensorwp}
g=\pi^{*}g_{B}+(f\circ \pi)^2 \sigma^{*}g_{F},
\end{equation}
is called warped product. We denote it by $B^{n}\times_{f} F^{m}$. where * is the pullback of canonical projections $\pi: B^{n}\times F^{m}\rightarrow B$, $\sigma: B^{n}\times F^{m}\rightarrow F$ which are projections on the first and second manifold, respectively. The function $f$ is called warping function, $B$ is called the base and $F$ the fiber.
\end{definition}

The authors in \cite{de2017gradient} studied gradient Ricci soliton on warped product and proved that either the warping function is constant or the potential function satisfies 
\begin{equation}\label{eqpotencialbase}
h=h_B\circ\pi, \quad h_B \in C^{\infty}(B).
\end{equation}

In this paper,  we will consider the harmonic map as a real function $u: M \rightarrow \mathbb{R}$.
Our first result characterize locally the harmonic map $u$ by means of the potential function $h$. More precisely, we obtain the following
\begin{proposition}\label{prop1}
Let $(B^n\times_fF^m,g,h,u,\lambda)$ be a GRHS with non-constant harmonic map $u$, then in a neighborhood $\mathcal{V}$ of a point $(p,q) \in B^{n}\times F^{m}$ the harmonic map $u$ can be represented as $u=u_B\circ \pi$ or $u=u_F\circ \sigma$ if, and only if, $h=h_B\circ\pi$ .
\end{proposition}

Motivated by the above proposition, we will consider the study of harmonic map in two cases, $u=u_B\circ \pi$ or $u=u_F\circ \sigma$. As a result we obtained the necessary and sufficient conditions for existence of the GRHS on warped product.

\begin{theorem}\label{theoequationfund}
$(B^n\times_fF^m,g,h,u,\lambda)$ is a GRHS if, and only if  $f,h,u,\lambda$ verify:
\begin{itemize}
\item (a) If $u=u_B\circ \pi,$ then
\end{itemize}
\begin{flalign}\label{grhsinbase}
\begin{cases}
&Ric_{g_B}-\frac{m}{f}Hess_{g_B}f + Hess_{g_B} h_B - \theta \nabla_{g_B} u_B \otimes \nabla_{g_B} u_B=\lambda g_{B},\\
&\Delta_{\omega}u_B=0\ \mbox{in}\ B,\\
\end{cases}
\end{flalign}

\begin{eqnarray*}
&F \ \mbox{is Einstein with}\ Ric_{g_F}=\mu g_F,
\end{eqnarray*}

\begin{eqnarray}\label{constmu}
&f\Delta_{g_{B}}f+(m-1)|\nabla_{g_B} f|^2+\lambda f^2+f\nabla_{g_B} f(h_B)=\mu.
\end{eqnarray} 

\begin{itemize}
\item (b) If $u=u_F\circ \sigma,$ then
\end{itemize}
\begin{flalign}\label{eqeinsharm}
&Ric_{g_B}-\frac{m}{f}Hess_{g_B} f + Hess_{g_B} h_B =\lambda g_{B},
\end{flalign}
\begin{flalign}
F\ \mbox{is Einstein-harmonic with}\ \begin{cases}Ric_{g_F}-\theta \nabla_{g_F} u_F \otimes \nabla_{g_F} u_F&=\mu g_F,\\
\Delta_{g_F}u_F=0,\end{cases}\nonumber
\end{flalign}
\begin{flalign}
f\Delta_{g_{B}}f+(m-1)|\nabla_{g_B} f|^2+\lambda f^2+f\nabla_{g_B} f(h_B)&=\mu.\label{constrho}
\end{flalign}
Where $Ric_{g_B}$, $Ric_{g_F}$ are the Ricci tensors of $B$, $F$ respectively, and $\Delta_{\omega}=\Delta-<\nabla, \nabla \omega>$ and $\omega=h-m\nabla log(f).$
\end{theorem}
\begin{observation}Note that the Theorem \ref{theoequationfund} is a generalization of corollary 3 obtained in \cite{10.2307/1194287}.
\end{observation}
In sequel, we obtain some triviality results for the warping function, potential function and harmonic map by means of the maximum principle.
\begin{theorem}\label{theotrivi}
Let $(B^n\times_fF^m,g,h,u,\lambda)$ a GHRS with non-constant $u$, then:
\begin{enumerate}
\item If $u_B,\ u_F$ reaches the maximum or minimum in $B,\ F$ respectively, then $u=u_B \circ \pi$ or $u=u_F \circ \sigma$  is a constant function, therefore $(B^n\times_fF^m,g,h,u,\lambda)$ is a gradient Ricci soliton,
\item If $\lambda \geq 0$, $h_B$ reaches the maximum or minimum in $B$ and $ \frac{m \Delta_{g_B} f}{f} \geq scal_{g_B}$, then $h=h_B \circ \pi$ is a constant function, therefore $(B^n\times_fF^m,g,h,u,\lambda)$ is a harmonic Einstein manifold,
\item If $f_B$ reaches the maximum and $\lambda \leq \frac{\mu}{f^2}$  in $B$, then $f=f_B \circ \pi$ is a constant function, therefore $(B^n\times_fF^m,g,h,u,\lambda)$ is a semi-Riemannian product. (The same occurs if $f_B$ reaches the minimum and $\lambda \geq \frac{\mu}{f^2}$).
\end{enumerate}
Where, $scal_B$ is scalar curvature of $B$.
\end{theorem}

In order to obtain infinitely many examples of geodesically complete steady GRHS we will use as a tool warped products semi-Riemannian metric conformal to Euclidean space as used  in \cite{MR3912818} in which the authors to expose a family of geodesically complete steady gradient Yamabe soliton. 

In the case where $u=u_B \circ \pi$, we consider the warped product semi-Riemannian GRHS with the base conformal to an $n -$dimensional semi-Euclidean space, invariant under the action of an $(n-1)$-dimensional translation group, that is, we take $(\mathbb{R}^n,\varphi^{-2} g_0)\times_{f}F^{m}$, where $g_0$ is the canonical semi-Riemannian metric and $\varphi$ is the conformal factor. More precisely, consider the semi-Riemannian metric $(g_0)_{ij}=\epsilon_i \delta_{ij}$ in local coordinates $x = (x_1, . . . , x_n)$ of $\mathbb{R}^n$, where $\epsilon_i = \pm 1.$ For an arbitrary choice of non zero vector $\alpha = (\alpha_1, . . . , \alpha_n)$ we define the function $\xi:\mathbb{R}^n \rightarrow \mathbb{R}$ by
\begin{eqnarray*}
\xi(x_1, . . . , x_n) = \sum_{i=1}^{n}\alpha_ix_i.
\end{eqnarray*}

In the next we obtain conditions for the functions $h\circ\xi,\ f\circ\xi,\ u\circ\xi$ and $\varphi\circ\xi$ to satisfy (\ref{definitiongrhs}) with $M=\mathbb{R}^n\times F^m$ and metric tensor $$g=\frac{1}{\varphi^2}g_0+f^2g_F.$$ 

We will omit the composition of these functions above with the translation functions. Moreover, the next result gives us a necessary and sufficient condition for the existence of GRHS warped product with base conformal to a semi-Euclidean space invariant by  action of an $(n-1)$-dimensional translation group when $u =u_B \circ \pi$.
\begin{theorem}\label{theoinvarinbase}
The $(\mathbb{R}^n\times_fF^m,g,h,u=u_B\circ\pi,\lambda)$ is a GRHS with harmonic map non-constant and $f=f\circ \xi,\ h=h\circ \xi,\ \varphi\circ\xi,\ u=u\circ\xi$ defined in $(\mathbb{R}^n,\varphi^{-2}g_0)$ furnished with the metric tensor $g=\varphi^{-2}g_0+f^2g_F$ if, and only if, the functions verify the system below:
\begin{equation}
(n-2)\frac{\varphi''}{\varphi}-m\frac{f''}{f}-2m\frac{\varphi'}{\varphi}\frac{f'}{f}+h''+2\frac{\varphi'}{\varphi}h'-\theta(u')^2=0,\label{eqinvbase1}\\
\end{equation}
\begin{equation}
\left[\frac{\varphi''}{\varphi}-(n-1)\left(\frac{\varphi'}{\varphi}\right)^2+m\frac{\varphi'}{\varphi}\frac{f'}{f}-\frac{\varphi'}{\varphi}h'\right]||\alpha||^2=\frac{\lambda}{\varphi^2},\label{eqinvbase2}\\
\end{equation}
\begin{equation}
\left[\frac{f''}{f}-(n-2)\frac{\varphi'}{\varphi} \frac{f'}{f}+(m-1)\left(\frac{f'}{f}\right)^2-\frac{f'}{f}h'\right]||\alpha||^2=\frac{\mu}{f^2\varphi^2}-\frac{\lambda}{\varphi^2},\label{eqinvbase3}\\
\end{equation}
\begin{equation}
\left[u''-(n-2)\frac{\varphi'}{\varphi}u'+mu'\frac{f'}{f}-u'h'\right]||\alpha||^2=0.\label{eqinvbase4}
\end{equation}
\end{theorem}

\begin{cor}
Under the hypotheses of Theorem \ref{theoinvarinbase} if $||\alpha||^2=0$ then, the GRHS\\ $(\mathbb{R}^n\times_f\mathbb{R}^m,g,h,u,\lambda)$ is steady and $F^m$ is Ricci flat with metric tensor $g=\varphi^{-2}g_0+f^2g_F$, that is, $\lambda=0$ and $\mu=0$.
\end{cor}

\begin{example}\label{examplecomplete}
In Theorem \ref{theoinvarinbase} consider $||\alpha||^2=0$ and  the Lorentzian space $(\mathbb{R}^{n},\overline{g})$ with coordinates $(x_{1},\dots,x_{n})$, signature $\varepsilon_{1}=-1$, $\varepsilon_{i}=1$, $\forall$ $i\geq2$, and fiber $(\mathbb{R}^{m},g_{0})$ where $g_{0}$ is the semi-Euclidean metric. Let $\xi=x_{1}+x_{2}$ and choose $k\in\mathbb{R}\setminus\{0\}$. Then
\begin{equation}
f(\xi)=e^{k\xi},\hspace{0,3cm}\varphi(\xi)=e^{k\xi},\hspace{0,3cm} u(\xi)=k\xi\nonumber,
\end{equation}
\begin{equation}
h(\xi)=\frac{k}{2}(2 - n + 3m + \theta)\xi  - \frac{e^{-2k \xi} k_1}{2 k} + k_2\nonumber
\end{equation}
defines a family of geodesically complete steady GRHS on $(\mathbb{R}^{n},\varphi^{-2}\overline{g})\times_{f}(\mathbb{R}^{m},g_{0})$ with potential
function $h\circ\pi$, harmonic map $u\circ\pi$ and warping function $f\circ\pi$.
\end{example}
\begin{example}
In Theorem \ref{theoinvarinbase} consider $||\alpha||^2=1$, $m,n=1$ and $k^2> 1$ then the functions 
\begin{equation}
h(\xi)=(m-k(n - 2))\ln( f(\xi)),\hspace{0,3cm}\varphi(\xi)=f(\xi)^k,\hspace{0,3cm} u(\xi)=k\xi\nonumber,
\end{equation}
\begin{equation}
f(\xi)=e^{\sqrt{\frac{\theta k^2}{(-k^2(n-2)-m)}} \xi},\nonumber
\end{equation}
provide a steady GRHS on $(\mathbb{R}^{n},\varphi^{-2}\delta)\times_{f}(\mathbb{R}^{m},\delta)$ defined for all  $\xi \in \mathbb{R}$. Were $\delta$ is standard metric Euclidean.
\end{example}

When $u=u_F \circ \sigma$ we consider warped product GRHS with the base and fiber conformal to an $n -$dimensional and $m-$dimensional semi-Euclidean spaces, invariants under the action of an $(n-1)$-dimensional and $(m-1)$-dimensional translation group respectively, that is $(\mathbb{R}^n,\varphi^{-2} g_0)\times_{f}(\mathbb{R}^m,\tau^{-2} g_0)$ where $\tau$ is conformal factor of the fiber. In this way, similarly we define $\zeta: \mathbb{R}^m \rightarrow \mathbb{R}$ by $\zeta(x_{n+1},...,x_{n+m})=\beta_{n+1}x_{n+1}+...+\beta_{n+m}x_{n+m}$ for an arbitrary choice of non zero vector $\beta = (\beta_{n+1}, . . . , \beta_{n+m})$ and $y=(x_{n+1},...,x_{n+m})$ in $\mathbb{R}^m$ with semi-Riemannian metric $g_0$.

In the next result, we obtain conditions for the functions  $h\circ\xi,\ f\circ\xi,\ u\circ\zeta,\ \varphi\circ\xi,\ \tau\circ\zeta$  to satisfy (\ref{definitiongrhs}) with $M=\mathbb{R}^n\times \mathbb{R}^m$ and metric tensor $$g=\frac{1}{\varphi^2}g_0+f^2\frac{1}{\tau^2}g_0.$$
\begin{theorem}\label{theoinvarinfiber}
The $(\mathbb{R}^n\times_f\mathbb{R}^m,g,h,u=u_F\circ\sigma,\lambda)$ is a GRHS with harmonic map non-constant and $f=f\circ \xi,\ h=h\circ \xi,\ \varphi\circ\xi,\ u=u\circ\zeta$ defined in $(\mathbb{R}^n,\varphi^{-2}g_0)$ and $(\mathbb{R}^m,\tau^{-2}g_0)$ furnished with the metric tensor $g=\varphi^{-2}g_0+f^2\tau^{-2}g_0$ if, and only if, the functions verify the system below:
\begin{equation}
(n-2)\frac{\varphi''}{\varphi}-m\frac{f''}{f}-2m\frac{\varphi'}{\varphi}\frac{f'}{f}+h''+2\frac{\varphi'}{\varphi}h'=0,\label{eqinvfiber1}
\end{equation}
\begin{equation}
\left[\frac{\varphi''}{\varphi}-(n-1)\left(\frac{\varphi'}{\varphi}\right)^2+m\frac{\varphi'}{\varphi}\frac{f'}{f}-\frac{\varphi'}{\varphi}h'\right]||\alpha||^2=\frac{\lambda}{\varphi^2},\label{eqinvfiber2}
\end{equation}
\begin{equation}
\begin{split}
\left[f''\varphi^2f-(n-2)\varphi'\varphi ff'+(m-1)(f')^2\varphi^2-f'f\varphi^2h'\right]||\alpha||^2+\lambda f^2=\\
=\left[\frac{\tau''}{\tau}-(m-1)\left(\frac{\tau'}{\tau}\right)^2\right]||\beta||^2,\label{eqinvfiber3}\\
\end{split}
\end{equation}
\begin{equation}
(m-2)\frac{\tau''}{\tau}-\theta(u')^2=0,\label{eqinvfiber4}
\end{equation}
\begin{equation}
\left[\tau^2u''-(m-2)\tau\tau'u'\right]||\beta||^2=0.\label{eqinvfiber5}
\end{equation}
\end{theorem}

\begin{cor}
Under the hypotheses of Theorem \ref{theoinvarinfiber} if $||\alpha||^2=0,$ then the GRHS $(\mathbb{R}^n\times_f\mathbb{R}^m,g,h,u=u_F\circ\sigma,\lambda)$ is steady with metric tensor $g=\varphi^{-2}g_0+f^2\tau^{-2}g_0$, that is, $\lambda=0$.
\end{cor}
In the next results we describe all the solutions of Theorem \ref{theoinvarinfiber} $\mathbb{R}^n\times_f\mathbb{R}^m$  when $g=\varphi^{-2}g_0+f^2\tau^{-2}g_0$ is a steady GRHS, i.e., $\lambda = 0,$ and $(\mathbb{R}^m,\tau^{-2}g_0)$ is a steady Einstein-harmonic manifold, i.e., $\rho=0.$
\begin{theorem}\label{theoclassiinvarinfiber}
Under considerations of Theorem \ref{theoinvarinfiber} the $\mathbb{R}^n\times_f\mathbb{R}^m$ is a steady GRHS with $u$ non-constant and $F$ Einstein-harmonic steady(i.e $\rho=0$), $n>1, m \geq 3,$ furnished with the metric tensor $g=\varphi^{-2}g_0+f^2\tau^{-2}g_0$ if, and only if, the functions $\varphi(\xi), f(\xi), h(\xi)$, $\tau(\zeta), u(\zeta)$ verify
\begin{enumerate}
\item If $||\alpha||^2=0, \quad ||\beta||^2=0$
\begin{equation}\label{potentialfunctionalphanull}
h(\xi) =\int{\left(\varphi^{-2}\left[\int{(m\frac{f''}{f}\varphi^2 + 2m\varphi\varphi'\frac{f'}{f}-(n-2)\varphi\varphi'')}d\xi + c_1\right]\right)}d\xi+c_2,
\end{equation}
\begin{equation}\label{harmonicmapfibranulo}
u(\zeta) =\int{\pm\sqrt{\frac{(m-2)}{\theta}\frac{\tau''}{\tau}(\zeta)}}d\zeta+c_3,
\end{equation}
where $c_1,c_2>0,c_3 \in \mathbb{R}.$
\item If $||\alpha||^2=0, \quad ||\beta||^2=1$, the potential function $h(\xi)$ is defined for \eqref{potentialfunctionalphanull}, and

\begin{equation}\label{eqfatorconformefiber}
\tau(\zeta) = c_2(c_1+(m-2)\zeta)^{\frac{1}{2-m}},
\end{equation}

\begin{equation}\label{equinfiberinvariant}
u(\zeta) = c_2\pm\frac{\sqrt{2 - 3m+m^2} ((-2 + m)\zeta+c_1)\log{((-2 + m)\zeta+c_1})}{(-2 + m) \sqrt{\theta((-2+m)\zeta+ c_1)^2}}.
\end{equation}

These solutions are defined on the half space defined by $\zeta>-\frac{c_1}{(m-2)}$.

\item If $||\alpha||^2=1, \quad ||\beta||^2=0$
\begin{equation}\label{eqbasenulofibernonulo}
\begin{split}
u(\zeta) =\int{\pm\sqrt{\frac{(m-2)}{\theta}\frac{\tau''}{\tau}(\zeta)}}d\zeta+c_3,
\end{split}
\end{equation}
\begin{equation}
\begin{cases}\label{functionwithzconst} 
\varphi_{\pm}(\xi)=\frac{c_4}{(N_{\pm}\xi+b)^{\frac{k}{N_{\pm}}}},\\ 
f_{\pm}(\xi)=\frac{c_5}{(N_{\pm}\xi+b)^{\frac{1}{N_{\pm}}}},\\ 
h_{\pm}(\xi)=-\frac{m-(n-2)k+N_{\pm}}{N_{\pm}}\ln(N_{\pm}\xi+b),
\end{cases}
\end{equation}
\end{enumerate}
when $h'(\xi)+[(n-2)k-m](f'/f)(\xi)=z(\xi)(f'/f)(\xi),$ where $z(\xi)$ is a nonzero constant, $c_4,c_5>0,k>0,b$ and $N_{\pm}$ be constant with $N_{\pm}=-k\pm\sqrt{m+k^2(n-1)}$. When $z(\xi)$ is a  smooth non constant function, the function $\varphi(\xi),f(\xi),h(\xi)$ are obtained by integrating the system

\begin{equation}\label{functionwithzfunc}
\begin{cases} 
\frac{\varphi'}{\varphi}(\xi)=k\frac{f'}{f}(\xi),\\ 
\frac{f'}{f}(\xi)=\psi(\xi),\\ 
h'(\xi) = [z(\xi)+m-k(n-2)]\psi(\xi),
\end{cases}
\end{equation}

\begin{equation}\label{functionwithzfunc1}
\begin{cases} 
\psi(\xi)=c_6(z(\xi)+k-\sqrt{m+k^2(n-1)})^{\frac{a-1}{2}}(z+k+\sqrt{m+k^2(n-1)})^{\frac{-a+1}{2}},\\
z'(\xi)=-c_6(z(\xi)+k-\sqrt{m+k^2(n-1)})^{\frac{a+1}{2}}(z + k + \sqrt{m + k^2(n-1)})^{\frac{-a-1}{2}},
\end{cases}
\end{equation}
where $a,k,c_6 \ \in \mathbb{R}$, with $k,c_6\geq 0$ and $a = \frac{k}{\sqrt{m + k^2(n-1)}}.$

\begin{enumerate}
\setcounter{enumi}{3}
\item If $||\alpha||^2=1, \quad ||\beta||^2=1$
\end{enumerate}
The functions $\varphi(\xi),f(\xi),h(\xi)$ are defined by the system \eqref{functionwithzconst} or \eqref{functionwithzfunc}, \eqref{functionwithzfunc1}. While the functions $u(\zeta)$ and $\tau(\zeta)$ are determined by equations \eqref{eqfatorconformefiber}, \eqref{equinfiberinvariant}.
\end{theorem}
\begin{example}\label{examplecomplefibir}
In the Theorem \ref{theoclassiinvarinfiber} Case(1) consider the Lorentzian spaces $(\mathbb{R}^{n},g_1)$, $(\mathbb{R}^{m},g_2)$, with coordinates $(x_{1},\dots,x_{n})$, $(x_{n+1},\dots,x_{n+m})$ where $g_1=-dx_{1}^{2}+\sum_{i=2}^{n}dx_{i}^{2}$ and $g_2=-dx_{n+1}^{2}+\sum_{j=n+2}^{n+m}dx_{j}^{2}$. Let $\xi=x_{1}+x_{2}$ and $\zeta=x_{n+1}+x_{n+2}$. Then
\begin{equation*}
f(\xi)=ke^{A\xi},\hspace{0,3cm}\varphi(\xi)=ke^{A\xi},\hspace{0,3cm} \tau(\zeta)=\zeta^2+1,\quad A\neq 0,\  k>0
\end{equation*}
\begin{equation*}
h(\xi)=(2 - n + 3m + \theta)\frac{A}{2}\xi  - \frac{e^{-2A \xi} c_7}{2Ak^2} + c_8,\quad  c_7,c_8\in\mathbb{R},
\end{equation*}
\begin{eqnarray*}
\begin{split}
u(\zeta)&= -\frac{\sqrt{2} \sqrt{-2+m} \sqrt{1+\zeta^2} \sinh^{-1}(\zeta)}{\sqrt{\left(1+\zeta^2\right) \theta }}+c_9,\quad c_9\in\mathbb{R},
\end{split}
\end{eqnarray*}
defines a family of geodesically complete steady GRHS on $(\mathbb{R}^{n},\varphi^{-2}g_1)\times_{f}(\mathbb{R}^{m},\tau^{-2}g_2)$ with potential function $h\circ\pi$, harmonic map $u\circ\sigma$ and warping function $f$(see section \ref{proofs}).
\end{example}
\begin{example}
In Theorem \ref{theoclassiinvarinfiber} Case(3) consider $k=1$, $z(\xi)$ a constant nonzero, the base $(\mathbb{R}^{2},g_0)$, the fiber Lorentzian spaces $(\mathbb{R}^{3},g)$ with coordinates $(y_{1},y_{2},y_{3})$ and $\zeta=y_{1}+y_{2}$. Than we obtain
\begin{equation*}
\begin{cases}
\varphi(\xi)=\frac{c_5}{\xi + b},\\ 
f(\xi)=\frac{c_4}{\xi + b},\\ 
h(\xi)=-4\ln(\xi+b),
\end{cases}
\end{equation*}
\begin{equation*}
\begin{cases}
\tau(\zeta)=e^{c_4\zeta},\\ 
u(\zeta)=\frac{c_4}{\sqrt{\theta}}\zeta.
\end{cases}
\end{equation*}
Which describes a family of steady GRHS defined in half-space $\xi>-b$ with metric tensor $g=\varphi^{-2}g_0+f^2\tau^{-2}g.$
\end{example}

\section{Proof of the main results}
\label{proofs}
\begin{myproof}{Proposition}{\ref{prop1}}
Let $\mathcal{L}(B)$, $\mathcal{L}(F)$ the spaces of lifts of vector fields on $B$ and $F$ to $B\times F$, respectively. Consider $X,Y\in\mathcal{L}(B)$ and $V,W\in\mathcal{L}(F)$. We have by the well known formula of Ricci curvature and hessian tensor on warped product $B^n\times_fF^m$ \cite{o1983semi} that
\begin{equation}\label{eqricciwp}
\begin{split}
Ric(X,Y)&=Ric_{g_B}(X,Y)-\frac{m}{f}Hess_{g_B} f (X,Y)\\
Ric(X,V)&=0\\
Ric(V,W)&= Ric_{g_F}(V,W)-\left(f\Delta_{B}f+ (m - 1) g_B(\nabla f,\nabla f) \right)g_F(V,W),
\end{split}
\end{equation}
 for all  $X,Y \in \mathcal{L}(B)$ and $V,W \in \mathcal{L}(F)$.

Suppose $u$ can be represented as $u=u_B\circ \pi$ or $u=u_F\circ \sigma$, then using \eqref{definitiongrhs}
\begin{equation*}
Hess h(X,V) = X (V (h))-\nabla_{X}V(h) = \theta\nabla u(X)\nabla u(V)=0,\quad \forall\ X \in \mathcal{L}(B),\ \forall\ V \in \mathcal{L}(F),
\end{equation*}
follow that by \cite{de2017gradient} that $h=h_B\circ\pi.$ On the other hand assume (\ref{eqpotencialbase}), then the equation (\ref{definitiongrhs}) applied to $X, V$ is
\begin{equation}\label{equbaseorfiber}
0=\lambda g(X,V)=\theta \nabla u(X)\nabla u(V),
\end{equation}	
because the metric tensor (\ref{metrictensorwp}) applied to $X,V$ is null. 

By hypothesis the harmonic map is not constant, it follows that there is a field\\ $L = X+V \in \mathcal(B^n\times F^m)$, $(p,q) \in B^n\times F^m$ and a neighborhood $\mathcal{V}$ such that \\$\nabla u(X+V)\nabla u(X+V)\neq 0$ in $\mathcal{V}$, than 
\begin{equation}
\nabla u(X)^2+2\nabla u(X)\nabla u(Y)+\nabla u(Y)^2\neq 0, 
\end{equation}
of (\ref{equbaseorfiber}) we have to
\begin{equation}
\nabla u(X)^2+\nabla u(Y)^2\neq 0.
\end{equation}
Therefore, $\nabla u(X)\neq 0$ or $\nabla u(Y)\neq 0$ in $\mathcal{V}$.
\end{myproof}
\begin{myproof}{Theorem}{\ref{theoequationfund}}

Assuming $u=u_B\circ\pi,$ by (\ref{eqricciwp}) and (\ref{definitiongrhs}) aplied in $X,Y \in \mathcal{L}(B)$ we get  \eqref{grhsinbase}. Similarly, applying (\ref{definitiongrhs}) in $V,W \in \mathcal{L}(F)$ combined with (\ref{eqricciwp}) we have to
\begin{eqnarray}\label{eqfundinfiber}
\lambda f^2g_F(V,W)-Hess h(V,W)&=&Ric_{g_F}(V,W)\nonumber\\
&-&\left(f\Delta_{B}f+ (m - 1) g_B(\nabla f,\nabla f) \right)g_F(V,W).
\end{eqnarray}
For all $V,W \in \mathcal{L}(F)$. Now, we obtain the Hessian expression $Hess h(V,W)$
\begin{eqnarray}
Hess h(V,W)&=&V(W(h))-(\nabla_{V}W)(h)\nonumber\\
&=&-(^F\nabla_{V}W-fg_F(V,W)\nabla_{g_B} f)(h)\nonumber\\
&=&fg_F(V,W)\nabla_{g_B} f(h)\label{hessinfiber}
\end{eqnarray}
Substituting (\ref{hessinfiber}) in (\ref{eqfundinfiber}) it follows that $Ric_{g_F}=\mu g_F$ were
\begin{equation}
\mu= f\Delta_{g_{B}}f+(m-1)|\nabla_{g_B} f|^2+\lambda f^2+f\nabla_{g_B} f(h).
\end{equation}
Therefore, $F$ is Einstein manifold. The reciprocal is immediate.

When $u=u_F\circ\sigma$ the proof is similar.
\end{myproof}

\begin{myproof}{Theorem}{\ref{theotrivi}}
\textbf{Case 1)} Using Proposition 35 of \cite{o1983semi} we have that
\begin{equation}\label{laplaubase}
\Delta u = \left[\Delta_{g_B} u_B+\frac{m}{f} g_B(\nabla_{g_B} u_B, \nabla_{g_B} f)\right] \circ \pi
\end{equation}
when $u=u_B\circ\pi$, and
\begin{equation}\label{laplaufiber}
\Delta u = \frac{\Delta_{g_F}u_F}{f^2} \circ \sigma\nonumber
\end{equation}
when $u=u_F\circ\sigma$.
Thus, the second equation in definition (\ref{definitiongrhs}) becomes
\begin{equation}\label{laplacianuinbase}
 \left[\Delta_{g_B} u_B+\frac{m}{f} g_B(\nabla_{g_B} u_B, \nabla_{g_B} f))\right] \circ \pi= g_B(\nabla u_B, \nabla h_B) \circ \pi,
\end{equation}
if $u=u_B\circ\pi$ or
\begin{equation}\label{laplacianuinfiber}
 \frac{\Delta_{g_F}u_F}{f^2} \circ \sigma= 0,  
\end{equation}
if $u=u_F\circ\sigma$.

Then applying the maximum principle(or minimum) (see \cite{gilbarg2015elliptic}) to (\ref{laplacianuinfiber}) we obtain that $u_F$ is constant, thus $u=u_F\circ\sigma$ is constant and $B^n\times_fF^m$ is gradient Ricci soliton.

Denoting $\Delta_{\omega}:=\Delta - g(\nabla\omega,\nabla\cdot)$ the equation (\ref{laplacianuinbase}) becomes
\begin{equation}
\Delta_{\omega}u_B=0\nonumber,
\end{equation}
where $\omega=h_B-mlog(f)$. It follows from the maximum principle(or minimum) that $u_B$ is constant and therefore $u=u_B\circ\pi$ is constant. 

\textbf{Case 2)}  Using (\ref{eqricciwp}) for $X,Y \in \mathcal{L}(B)$ in (\ref{definitiongrhs}) follow that
\begin{equation}
\lambda g_B(X,Y)- Hess_{g_B} h(X,Y)+ \theta \nabla_{g_B}u(X)\nabla_{g_B}u(Y)= Ric_{g_B}(X,Y)-\frac{m}{f}Hess_{g_B} f (X,Y),\nonumber
\end{equation}
taking trace in relation to metric of $B$ on both sides gives
\begin{equation}
\Delta_{g_B}h_B=n\lambda+\theta||d\pi(\nabla u)||^2-scal_{g_B} + \frac{m \Delta_{g_B}f}{f}.\nonumber
\end{equation}
Therefore, using the hypotheses follow that $\Delta_{g_B}h_B\geq 0$ and maximum principle we have that $h_B$ is constant, consequently $h=h_B\circ\pi$ is constant function.
\begin{observation}Note that in the demonstration of case 2 either $u=u_B\circ\pi$ or $u=_F\circ\sigma$ the inequality is valid.
\end{observation}
\textbf{Case 3)} Consider the elliptic operator of second order given by
\begin{equation*}
\Xi(\cdot)=\Delta(\cdot)-\nabla h (\cdot)+\frac{(m-1)}{f}\nabla f(\cdot).
\end{equation*}
By the (\ref{constmu}) and (\ref{constrho}) of Theorem \ref{theoequationfund} follow that 
\begin{equation}
f\Delta_{g_{B}}f+(m-1)|\nabla_{g_B} f|^2+\lambda f^2+f\nabla_{g_B} f(h)=\eta
\end{equation}
or equivalently
\begin{equation}
\Xi(f)=\frac{\eta-\lambda f^2}{f}.
\end{equation}
Therefore, by hypothesis and strong maximum(or minimum) principle follow the result.
\end{myproof}

\begin{myproof}{Theorem}{\ref{theoinvarinbase}}
The Theorem \ref{theoequationfund} gives us necessary and sufficient condition to $B^n\times_fF^m$ be a GRHS, thus we will use it combined with invariant solution technique to obtain equations (\ref{eqinvbase1}), (\ref{eqinvbase2}), (\ref{eqinvbase3}) and (\ref{eqinvbase4}).

First, for an arbitrary choice of a non zero vector $\alpha=(\alpha_{1},\dots,\alpha_{n})$, consider  $\xi:\mathbb{R}^{n}\rightarrow\mathbb{R}$ given by $\xi(x_{1},\dots,x_{n})=\alpha_{1}x_{1}+\dots+\alpha_{n}x_{n}$. Since we are assuming that $\varphi(\xi)$, $h(\xi)$, $u(\xi)$ and $f(\xi)$ are functions of $\xi$, then we have
	\begin{equation}\label{invariante}
	\begin{aligned}
	\varphi_{,x_{i}}&=\varphi'\alpha_{i},\hspace{0.2cm} & f_{,x_{i}}&=f'\alpha_{i},\hspace{0.2cm} & h_{,x_{i}}&=h'\alpha_{i},\hspace{0.2cm} & u_{,x_{i}}&=u'\alpha_{i}\\[10pt]
	\varphi}_{,x_{i}x_{j}&=\varphi''\alpha_{i}\alpha_{j},\hspace{0.2cm} & f_{,x_{i}x_{j}}&=f''\alpha_{i}\alpha_{j},\hspace{0.2cm} &h_{,x_{i}x_{j}}&=h''\alpha_{i}\alpha_{j}\hspace{0.2cm} & u_{,x_{i}x_{j}}&=u''\alpha_{i}\alpha_{j}.
	\end{aligned}
	\end{equation}
	
Remember that we are $f=f_B\circ\pi,\ \varphi=\varphi_B\circ\pi,\ h=h_B\circ\pi$ and $u=u_B\circ\pi,\ $ where $B=(\mathbb{R}^n,\varphi^{-2}g_0)$.
	It is well known that for the conformal metric  $g_{B}=\varphi^{-2}g_0$, the Ricci curvature is given by \cite{besse2007einstein}:	
\begin{equation*}\label{Ricciinvar}Ric_{g_{B}}=\frac{1}{\varphi^{2}}\Big{\{}(n-2)\varphi Hess_{g_0}(\varphi)+[\varphi\Delta_{g_0}\varphi-(n-1)|\nabla_{g_0}\varphi|^{2}]g_0\Big{\}}.
\end{equation*}
Sice  $(Hess_{g_0}(\varphi))_{i,j}=\varphi''\alpha_i\alpha_j,\ \Delta_{g_0}\varphi=\varphi''||\alpha||^2,$ and $ |\nabla_{g_0}\varphi|^2=\varphi'||\alpha||^2 $ we have

\begin{equation}\label{ricciineqj}
(Ric_{g_B})_{i,j} =\frac{1}{\varphi}\{(n-2)\varphi''\alpha_i\alpha_j\} \ \ \ \forall i\neq j=1,...,n
\end{equation}
\begin{equation}\label{ricciieqj}
(Ric_{g_B})_{i,i} =\frac{1}{\varphi^2}\left\{(n-2)\varphi\varphi''(\alpha_i)^2+[\varphi\varphi''||\alpha||^2-(n-1)(\varphi')^2||\alpha||^2]\epsilon_i\right\} \ \ \ \forall i=1,...,n.
\end{equation}

 Computing the $Hess (h)$ relatively to $g_B$ we have
	\begin{equation*}
	(Hess_{g_B}(h))_{ij}=h_{,x_{i}x_{j}}-\sum_{k=1}^{n}\Gamma_{ij}^{k}h_{,x_{k}},
	\end{equation*}
	where the Christoffel symbol $\Gamma_{ij}^{k}$ for distinct $i,j,k$ are given by
	\begin{equation}\Gamma_{ij}^{k}=0,\ \Gamma_{ij}^{i}=-\frac{\varphi_{,x_{j}}}{\varphi},\ \Gamma_{ii}^{k}=\varepsilon_{i}\varepsilon_{k}\frac{\varphi_{,x_{k}}}{\varphi}\;\ \mbox{and}\;\ \Gamma_{ii}^{i}=-\frac{\varphi_{,x_{i}}}{\varphi}.\nonumber
	\end{equation}
	Therefore,
		\begin{eqnarray}\label{hessian}(Hess_{g_B}(h))_{ij}&=&h_{,x_{i}x_{j}}+\varphi^{-1}(\varphi_{,x_{i}}h_{,x_{j}}+\varphi_{,x_{j}}h_{,x_{i}})-\delta_{ij}\varepsilon_{i}\sum_{k}\varepsilon_{k}\varphi^{-1}\varphi_{,x_{k}}h_{,x_{k}}\nonumber\\
	&=&\alpha_{i}\alpha_{j}h''+(2\alpha_{i}\alpha_{j}-\delta_{ij}\varepsilon_{i}||\alpha||^{2})\varphi^{-1}\varphi'h'.
	\end{eqnarray}
And the Laplacian $\Delta_{g_B} f=\sum_{k}\varphi^{2}\varepsilon_{k}(Hess_{g_B}(f))_{kk}$ of $f$ is
\begin{equation}\label{laplacianfunctorsion}
\Delta_{g_B} f=||\alpha||^{2}\varphi^{2}(f''-(n-2)\varphi^{-1}\varphi'f').
\end{equation}

On the other hand, the expression of $\nabla(f)h$, $|\nabla f|^{2}$ and $(\nabla u \otimes \nabla u)_{ij}$ on conformal metric $g_{B}$ are given by
\begin{eqnarray}\label{gradientes}
\begin{split}
\nabla_{g_B}(f)h=\langle\nabla_{g_B} f,\nabla_{g_B} h\rangle=\varphi^{2}\sum_{k}\varepsilon_{k}f_{,x_{k}}h_{,x_{k}}=||\alpha||^{2}\varphi^{2}f'h',\\
|\nabla_{g_B} f|^{2}=\varphi^{2}\sum_{k}\varepsilon_{k}f_{,x_{k}}^{2}=||\alpha||^{2}\varphi^{2}(f')^{2},\quad (\nabla_{g_B}u\otimes\nabla_{g_B}u)_{ij}=u_{,x_{i}}u_{,x_{j}}=\alpha_i\alpha_j(u')^2.
\end{split}
\end{eqnarray}
Then substituting \eqref{ricciieqj}, \eqref{hessian} and \eqref{gradientes} for $i=j$ into \eqref{grhsinbase} we obtain \eqref{eqinvbase2}.
	
	Now, for $i\neq j$ we obtain by \eqref{ricciineqj} and \eqref{hessian} that 
	\begin{equation*}
	\alpha_{i}\alpha_{j}\left((n-2)\frac{\varphi''}{\varphi}-m\frac{f''}{f}-2m\frac{\varphi'}{\varphi}\frac{f'}{f}+h''+2\frac{\varphi'}{\varphi}h'-\theta(u')^2\right)=0.
	\end{equation*}
	
	If there exist $i,j$, $i\neq j$ such that $\alpha_{i}\alpha_{j}\neq
	0$, then we get
	\begin{equation*}(n-2)\frac{\varphi''}{\varphi}-m\frac{f''}{f}-2m\frac{\varphi'}{\varphi}\frac{f'}{f}+h''+2\frac{\varphi'}{\varphi}h'-\theta(u')^2=0
	\end{equation*}
which is equation \eqref{eqinvbase1}. 

Now, we need to consider the case $\alpha_i\alpha_j=0, \ \forall i\neq j$. For this consider $k_0$ fixed and $\alpha_{k_{0}}=1$, $\alpha_{k}=0$ for $k\neq k_{0}$. In this case, substituting \eqref{ricciieqj}, \eqref{hessian} and \eqref{gradientes} into \eqref{grhsinbase} we obtain the equation \eqref{eqinvbase2} for $i\neq k_{0}$, that is, $\alpha_{i}=0.$ And we obtain the \eqref{eqinvbase1} when $i=k_{0}$, that is, $\alpha_{k_{0}}=1$.

Substituting \eqref{laplacianfunctorsion}, \eqref{gradientes} in \eqref{constmu} we obtain \eqref{eqinvbase3}. Finally, consider the second equation of the \eqref{definitiongrhs} combined with  \eqref{laplaubase} we have
\begin{equation}\label{homogenbase}
\Delta u = \left[\Delta_{g_B} u_B+\frac{m}{f} g_B(\nabla u, \nabla f)\right] \circ \pi=g_B(\nabla_{g_B}u,\nabla_{g_B}h)\circ\pi.
\end{equation}
Using \eqref{laplacianfunctorsion}, \eqref{gradientes} in \eqref{homogenbase} we obtain \eqref{eqinvbase4}, and this completes the demonstration.
\end{myproof}

\begin{myproof}{Theorem}{\ref{theoinvarinfiber}}
Now, we will be using the technique of invariant solutions in both the base and the fiber, thus by replacing the equation \eqref{ricciineqj} and \eqref{gradientes} in \eqref{eqeinsharm} for $i\neq j$ we obtain similarly the proof of above theorem the equation \eqref{eqinvfiber1} and if $i=j$ it is sufficient to use \eqref{ricciieqj}, \eqref{gradientes} in \eqref{eqeinsharm} for obtain \eqref{eqinvfiber2}.

We know by Theorem \ref{theoequationfund} that $F$ is Einstein-harmonic manifold, this is 
\begin{equation}\label{eqeinsteinharmonic}
Ric_{g_F}-\theta\nabla_{g_F}u_{g_F}\otimes\nabla_{g_F}u_{g_F}=\rho g_F.
\end{equation} Where $\theta>0$ and
\begin{equation}\label{eqpho}
\Delta_{g_{B}}f+(m-1)|\nabla_{g_B} f|^2+\lambda f^2+f\nabla_{g_B} f(h)=\rho.
\end{equation}
For arbitrary choice of a non zero vector $\beta=(\beta_1,...,\beta_m)$ consider $\tau:\mathbb{R}^m\rightarrow(0,\infty),\zeta:\mathbb{R}^m\rightarrow\mathbb{R}$ the conformal factor of the fiber and invariant function respectively. Since we are assuming $u(\zeta)$  then we have
\begin{equation}\label{eqtensoru}
(\nabla_{g_F}u_{g_F}\otimes\nabla_{g_F}u_{g_F})_{ij}= u_{,y_i}u_{,y_j}\beta_i\beta_j\quad \quad \forall i,j=1,...,m.
\end{equation}
Using \eqref{gradientes} in \eqref{eqpho} we have
\begin{equation}\label{rhoconstinvbase}
\left[f''\varphi^2f-(n-2)\varphi'\varphi ff'+(m-1)(f')^2\varphi^2-f'f\varphi^2h'\right]||\alpha||^2+\lambda f^2=\rho.
\end{equation}
Substituting \eqref{ricciineqj}, \eqref{ricciieqj}, \eqref{eqtensoru} and \eqref{rhoconstinvbase} in \eqref{eqeinsteinharmonic} we obtain the equations \eqref{eqinvfiber3} and \eqref{eqinvfiber4} for $i=j$ or $i\neq j$. Finally, by Theorem \ref{theoequationfund} part (b) we have $\Delta_{g_F}u_F=0$, then using \eqref{gradientes} we have \eqref{eqinvfiber5} which ends the proof of the theorem. 
\end{myproof}
\begin{myproof}{Theorem}{\ref{theoclassiinvarinfiber}} case 1): Sice $\lambda=0$ and $||\alpha||^2=0=||\beta||^2$ we have by Theorem \ref{theoinvarinfiber}
\begin{equation*}
(n-2)\frac{\varphi''}{\varphi}-m\frac{f''}{f}-2m\frac{\varphi'}{\varphi}\frac{f'}{f}+h''+2\frac{\varphi'}{\varphi}h'=0,
\end{equation*}
\begin{equation*}
(m-2)\frac{\tau''}{\tau}-\theta(u')^2=0,
\end{equation*}
this way, it is easy to see that these equations are \eqref{potentialfunctionalphanull}, \eqref{harmonicmapfibranulo}.

Case 2): By hypothesis $||\beta||^2=1$, $||\alpha||^2=0$, $\lambda=0$ and Theorem \ref{theoinvarinfiber} we have 
\begin{equation}\label{eqtau}
\begin{split}
\frac{\tau''}{\tau}-(m-1)\left(\frac{\tau'}{\tau}\right)^2=0,
\end{split}
\end{equation}
\begin{equation}\label{eqtaueu}
(m-2)\frac{\tau''}{\tau}-\theta(u')^2=0,
\end{equation}
\begin{equation}\label{eqtaueu1}
\tau^2u''-(m-2)\tau\tau'u'=0.
\end{equation} 

Integrating \eqref{eqtau} and substituing in \eqref{eqtaueu} or \eqref{eqtaueu1}  we obtain \eqref{eqfatorconformefiber} and \eqref{equinfiberinvariant}.

Case 3): Sice $||\alpha||^2=1,\ ||\beta||^2=0,\ \lambda=0$ we have by equations \eqref{eqinvfiber1},\eqref{eqinvfiber2},\eqref{eqinvfiber3} and \eqref{eqinvfiber4}  of Theorem \ref{theoinvarinfiber}

\begin{equation}\label{sys01}
\begin{cases}
(n-2)\frac{\varphi''}{\varphi}-m\frac{f''}{f}-2m\frac{\varphi'}{\varphi}\frac{f'}{f}+h''+2\frac{\varphi'}{\varphi}h'=0,\\
\frac{\varphi''}{\varphi}-(n-1)\left(\frac{\varphi'}{\varphi}\right)^2+m\frac{\varphi'}{\varphi}\frac{f'}{f}-\frac{\varphi'}{\varphi}h'=0,\\
f''\varphi^2f-(n-2)\varphi'\varphi ff'+(m-1)(f')^2\varphi^2-f'f\varphi^2h'=0,\\
(m-2)\frac{\tau''}{\tau}-\theta(u')^2=0.
\end{cases}
\end{equation}

The first three equations of system involving $\varphi,f,h$ was solved by \cite{de2017gradient}. While from the last equation of system \eqref{sys01} we obtain \eqref{eqbasenulofibernonulo}.

Case 4): It follows immediately from previous cases.
\end{myproof}
\begin{myproof}{completeness of example}{\ref{examplecomplefibir}}

\begin{definition}\label{defi}A semi-Riemannian manifold for which every geodesic is defined on the entire real line is said to be \textit{geodesically complete}, or just
\textit{complete}.
\end{definition}

Given a curve $\gamma$ in $M\times_{f}F$, we can write $\gamma(s)=(\gamma_{B}(s),\gamma_{F}(s))$, where $\gamma_{B}=\pi\circ\gamma$ and $\gamma_{F}=\sigma\circ\gamma$. The following proposition guarantees a condition for curve $\gamma$ to be geodesic.

\begin{proposition}\label{geodesicas produto torcido}(\cite{o1983semi}) A curve $\gamma=(\gamma_{B}, \gamma_{F})$ in $B\times_{f}F$ is a geodesic if, and only if,
	\\[0.1pt]
	\begin{enumerate}
		\item $\gamma_{B}''=g_{F}(\gamma_{F}',\gamma_{F}')f\circ\gamma_{B}\nabla f$\quad \text{in B},\\[0.1pt]
		\item $\gamma_{F}''=\dfrac{-2}{f\circ\gamma_{B}}\dfrac{d(f\circ\gamma_{B})}{ds}\gamma_{F}'$ \quad \text{in F}.
	\end{enumerate}
\end{proposition}

Let $(\mathbb{R}^{n},g_1)$, $(\mathbb{R}^{m},g_2)$ be the standard semi-Euclidean spaces as in the Example \ref{examplecomplefibir}. Take $A\in\mathbb{R}\setminus\{0\}, k>0$, and consider the functions \begin{equation}\varphi(\xi)=ke^{A\xi}, \hspace{0,3cm}f(\xi)=ke^{A\xi},\hspace{0,3cm}h(\xi)=(2 - n + 3m + \theta)\frac{A}{2}\xi  - \frac{e^{-2A \xi} c_7}{2Ak^2} + c_8,\quad  c_7,c_8\in\mathbb{R},\nonumber
	\end{equation}
and, 
\begin{eqnarray*}
\begin{split}
\tau(\zeta)=\zeta^2+1,\quad u(\zeta)&=-\frac{\sqrt{2} \sqrt{-2+m} \sqrt{1+\zeta^2} \sinh^{-1}(\zeta)}{\sqrt{\left(1+\zeta^2\right) \theta }}+c_9,\quad c_9\in\mathbb{R}.
\end{split}
\end{eqnarray*}
Where the function $f, \varphi, h$ are defined in $(\mathbb{R}^{n},g_1)$, end $\tau, u$ are defined in $(\mathbb{R}^{m},g_2)$.

We have to $g_B:=\varphi^{-2}g_1=ke^{-2A\xi}g_1$, then the gradient $\nabla_{g_B}f$ is given by
	\begin{equation*}
	\nabla_{g_B}f=\sum_{r,s=1}^{n}g_B^{rs}f_{,x_{s}}\partial_{r}=\sum_{r,s=1}^{n}\varphi^{2}\varepsilon_{r}\delta_{rs}f'\alpha_{s}\partial_{s}=\sum_{s=1}^{n}k^3A\varepsilon_{s}\alpha_{s}e^{3A\xi}\partial_{s}.
	\end{equation*}
Since $\alpha_{1}=\alpha_{2}=1$, $\alpha_{i}=0$, for $i\geq3$, and $\varepsilon_{1}=-1$, $\varepsilon_{i}=1$, for $i\geq2$, we obtain
	\begin{equation*}
\nabla_{g_B}f=\big{(}-k^3Ae^{3A\xi},k^3Ae^{3A\xi}, 0,\dots, 0\big{)}.
\end{equation*}
Then, considering $\gamma_{B}(s)=(y_{1}(s),\dots, y_{n}(s))$ and $\gamma_{F}(s)=(y_{n+1},\dots, y_{n+p}(s))$ in Proposition \ref{geodesicas produto torcido}, we have
\begin{equation*}
\left\{\begin{array}{r@{\mskip\thickmuskip}l}
	y_{1}''(s)&= -k^4Ag_F(\gamma_F',\gamma_F',)e^{4A(y_{1}(s)+y_{2}(s))}, \hfill \text{(I)}\\ [10pt]
	y_{2}''(s)&= k^4Ag_F(\gamma_F',\gamma_F',)e^{4A(y_{1}(s)+y_{2}(s))}, \hfill \text{(II)}\\[10pt]
	y_{r}''(s)&= 0,\quad \text{for}\quad r\in\{3,\dots,n\},\hfill \text{(III)}\\[10pt]
	y_{n+l}''(s)&=-2A[y_{1}'(s)+y_{2}'(s)]y_{n+l}'(s),\quad \text{for}\quad l\in\{1,\dots,d\}.\qquad \text{(IV)}
\end{array} \right.
\end{equation*}

Adding the differential equations \text{(I)} and \text{(II)} we have $y_{1}''(s)+y_{2}''(s)=0$, then by integration 
\begin{equation}\label{diferential eq1}
y_{1}'(s)+y_{2}'(s)=c_{1}, \quad y_{1}(s)+y_{1}(s)=c_{1}s+c_{2}, \quad c_{1}, c_{2}\in \mathbb{R}.
\end{equation}

Substituting \eqref{diferential eq1} into \text{(IV)}, we obtain the second order linear ordinary differential equation
\begin{equation*}
y_{n+l}''(s)+2Ac_{1}y_{n+l}'(s)=0 \quad \text{for each} \quad l\in\{1,\dots,d\},
\end{equation*}
whose general solutions is
\begin{equation}\label{second order}
y_{n+l}(s) = \begin{cases} c_{3,l}+c_{4,l}s &\text{if } c_{1}=0 \\[10pt]
	c_{3,l}+c_{4,l}e^{-2Ac_{1}s} & \text{if } c_{1}\neq0 
\end{cases}
\end{equation}
where $c_{3,l},c_{4,l}\in\mathbb{R}$. Therefore, for each $l\in\{1,\dots,m\}$, the functions $y_{n+l}(s)$ are defined on the entire real line $\mathbb{R}$. Notice that the solutions of \text{(III)} are given by $y_{r}(s)=c_{5,r}+c_{6,r}s$, for $c_{5,r}, c_{6,r}\in\mathbb{R}$, whose domain is also the real line, it is only necessary to prove that the solutions of \text{(I)} and \text{(II)} are also defined in $\mathbb{R}$, for this, we will explain $g_F(\gamma_F',\gamma_F')$ in function of the s parameter.

Since $g_F=\tau^{-2}g_1$, $\beta_{n+1}=\beta_{n+2}=1$, $\beta_{n+j}=0$, for $j\geq 3$, and $\varepsilon_{n+1}=-1$, $\varepsilon_{n+j}=1$, for $j\geq2$, we obtain

\begin{eqnarray*}
g_F(\gamma_F',\gamma_F')&=&\tau(\zeta\circ\gamma_F)^{-2}g_2(\gamma_F',\gamma_F')=\left[-y_{n+1}'(s)^2+y_{n+2}'(s)^2+\cdots y_{n+m}'(s)^2\right]\\
&=&\frac{1}{((y_{n+1}(s)+y_{n+2}(s))^2+1)^2}\left[-y_{n+1}'(s)^2+y_{n+2}'(s)^2+\cdots y_{n+m}'(s)^2\right],
\end{eqnarray*}
replacing \eqref{second order} into \text{(I)}, we have
\begin{equation*}
y_{1}''(s)=-k^4A\left[(c_{7,1}+c_{7,2}e^{-2Ac_1s})^2+1\right]^{-2}\left[-c_{8,1}^{2}+c_{8,2}^{2}+\dots+c_{8,r}^{2}\right]e^{-4Ac_{1}s}e^{4A(y_{1}(s)+y_{2}(s))}
\end{equation*}
where $c_{7,1}, c_{7,2}, c_{8,1}, c_{8,2},\dots,c_{8,r}\in\mathbb{R}$.

Now, by \eqref{diferential eq1} we obtain that

\begin{eqnarray}
y_{1}''(s)&=&c_{9,1}\left[(c_{7,1}+c_{7,2}e^{-2Ac_1s})^2+1\right]^{-2}e^{-4Ac_{1}s}e^{4A(c_{1}s+c_{2})}\nonumber\\
&=&c_{9,1}\left[(c_{7,1}+c_{7,2}e^{-2Ac_1s})^2+1\right]^{-2}e^{4Ac_{2}}\nonumber\\
&=& c_{10,1}\left[(c_{7,1}+c_{7,2}e^{-2Ac_1s})^2+1\right]^{-2} \quad c_{9,1}, c_{10,1}\in\mathbb{R}.
\end{eqnarray}
Note that the field $y_1'(s)$ is differentiable at every point and have a limited derivative. Indeed
\begin{eqnarray*}
|y_{1}''(s)|=\left|c_{10,1}\left[(c_{7,1}+c_{7,2}e^{-2Ac_1s})^2+1\right]^{-2}\right|<\left|c_{10,1}\right|.
\end{eqnarray*}
Thus, the system defined by
\begin{equation*}
\begin{cases}
y_1'(s)=z_1(s),\nonumber\\
z_1'(s)= \left[(c_{7,1}+c_{7,2}e^{-2Ac_1s})^2+1\right]^{-2},
\end{cases}
\end{equation*}
has solutions whose domain is $\mathbb{R}$. Except for the signal, the same occurs for $y_{2}(s)$. Thus, all the geodesics $\gamma=(\gamma_{B},\gamma_{F})$ are defined for the entire real line, which means that $(\mathbb{R}^{n},\varphi^{-2}g_1)\times_{f}(\mathbb{R}^{m},\tau^{-2}g_{2})$ is geodesically complete.
\end{myproof}

\bibliographystyle{abbrv}

\end{document}